\documentclass[12pt]{amsart}
\usepackage{latexsym, amsmath}

\theoremstyle{plain}
  \newtheorem{Thm}{Theorem}[section] 
  \newtheorem{Lma}[Thm]{Lemma} 
  \newtheorem{Cor}[Thm]{Corollary} 
  \newtheorem{Prop}[Thm]{Proposition}

\theoremstyle{definition}
  \newtheorem{Def}[Thm]{Definition}

\theoremstyle{remark}
  \newtheorem{Rem}[Thm]{Remark}
  \newtheorem{Q}[Thm]{Question}

\newcommand{\tensor}{\otimes}
\newcommand{\brq}{^{[q]}}

\newcommand{\mlabel}[1]%
  {\mbox{}\marginpar{\raggedleft\hspace{0pt}{\rm\ttfamily#1}}\label{#1}}

\newcommand{\depth}{\operatorname{depth}}

\newcommand{\length}{\operatorname{\lambda}}

\newcommand{\Hom}[3]{\operatorname{Hom}_{#1}(#2,#3){}}

\newcommand{\fr}{r_F}
\newcommand{\frp}{r_F^+}
\newcommand{\frm}{r_F^-}

\newcommand{\smi}{s^+}
\newcommand{\spl}{s^-}
\newcommand{\sdim}{\operatorname{sdim}}
\newcommand{\fm}{{\mathfrak m}}
\newcommand{\fn}{{\mathfrak n}}

\newcommand{\ringR}{\text{$(R,\fm,k)$ }}

 \newcommand{\height}{{\rm ht}}
 \newcommand{\ann}{{\rm Ann}}

\newcommand{\Dim}{{\rm dim}}

\newcommand{\Soc}{{\rm Soc}}
\newcommand{\E}{E_R(k)}
\newcommand{\Spec}{{\rm Spec}}

\newcommand{\et}{0^{\ast}_{E_{R}}}

\newcommand{\cs}{{\widetilde \tau}}
\newcommand{\im}{{\rm Im}}

\newcommand{\al}{\alpha}
\newcommand{\p}{\mathcal{P}}
\newcommand{\q}{\mathcal{Q}}

\newcounter{hours}\newcounter{minutes}

\newcommand{\excise}[1]{}
\begin{document}

\title{\bf The structure of F-pure rings}
\author[I.~M.~Aberbach]{Ian M. Aberbach}
\author[F.~Enescu]{Florian Enescu}
\address{Department of Mathematics, University of Missouri, Columbia, MO 65211, USA}
\email{aberbach@math.missouri.edu}
\address{Department of Mathematics, University of Utah, Salt Lake City,
UT  84112, USA and The~Institute of Mathematics of the Romanian Academy, Bucharest, Romania}
\email{enescu@math.utah.edu}
\thanks{The first author was partially funded by the National Security Agency.}
\date{}

\maketitle 

\begin{center}
{\it Dedicated to Professor Melvin Hochster on the occasion of his sixtieth birthday}
\end{center}

\begin{abstract}
For a reduced $F$-finite ring $R$ of characteristic $p >0 $ and $q=p^e$ one can
write $R^{1/q} = R^{a_q} \oplus M_q$, where $M_q$ has no free direct
summands over $R$. We investigate the structure of $F$-finite, $F$-pure rings
$R$ by studying how the numbers $a_q$ grow with respect to
$q$. This growth is quantified by the splitting dimension and the splitting
ratios of $R$ which we study in detail. We also prove the existence of
a special prime ideal $\p(R)$ of $R$, called the splitting prime, 
that has the property that $R/\p(R)$ is strongly $F$-regular. We show that this ideal captures
significant information with regard to the $F$-purity of $R$.
\end{abstract}

\section{Introduction} Let $\ringR$ be a reduced, local ring of positive
characteristic $p>0$ and Krull dimension $d$. Throughout this paper, $q=p^e$
denotes a power of the characteristic, $R^{1/q}$ is the ring of $q$th roots of
elements in $R$,  and $R$ will usually be assumed to be $F$-finite, i.e. $R^{1/p}$ is module
finite over $R$. This implies that $R^{1/q}$ is module finite  over $R$ for all
$q$. Write $R^{1/q} = R^{a_q} \oplus M_q$, which is a direct  sum decomposition
of $R^{1/q}$ over $R$ such that $M_q$ has no free direct  summands. We would
like to investigate the size of the numbers $a_q$, as $q$ grows  to infinity,
by studying the splitting dimension of the ring $R$ (for short the $s$-dimension
of $R$, as it was called in~\cite{AL}). (For precise definitions, we refer the
reader to Section 2).  The rings for which some $a_q > 0$ (equivalently,
all $a_q >0$) are called {\sl $F$-pure}.  The purpose of this paper is
to demonstrate that there is  much more structure to be discerned in such
rings, and at least part of this structure is captured by the $s$-dimension
and the splitting ratios (defined below).

The behavior of the sequence $\{ a_q \} _q$ captures subtle information about
the structure of the ring $R$ with regard to the natural action of Frobenius on
it. This behavior is intimately connected to the notions of $F$-purity and
strong $F$-regularity. In fact, our work starts with and extends ideas
originating in papers of  Huneke and Leuschke (\cite{HL}), and Aberbach and
Leuschke (\cite{AL}) where the  relation between the strong $F$-regularity and
the $F$-signature of  Gorenstein rings (\cite{HL}) and Cohen-Macaulay 
rings (\cite{AL}) are analyzed. One can think of the
concept of strong $F$-regularity as a very strong form of $F$-purity, as the
definitions included below show. Our work shows that there is a spectrum
of $F$-pure behavior, ranging from ``merely $F$-pure,'' in the
case that $\sdim(R) =0$ (i.e., essentially $a_q =1$ for all $q$), all
the way up to strong $F$-regularity (i.e., $\sdim(R) = \dim(R)$). 
Moreover, the central
notions of the paper, the splitting dimension and the Frobenius splitting 
ratios, will quantify this behavior in a clear fashion.

An important outcome of our analysis is the discovery, in any $F$-pure ring $R$,  of a special prime ideal
$\p(R)$, which we will call the \emph{splitting prime} of $R$, whose dimension (i.e. $\Dim(R/\p)$)
is an upper bound for the splitting
dimension of $R$. A significant feature of this prime ideal $\p$ is that it
defines a strongly $F$-regular quotient $R/\p$. In the case of a Gorenstein
ring $R$ this ideal is the only $F$-stable prime of $R$ (the study of the
$F$-stable primes of ring $R$ was initiated in (\cite{E}).

Let $\ringR$ be an $F$-finite $F$-pure ring.  For simplicity in making an
initial statement, assume that the residue field of $R$ is perfect.
 We let $k$ be the largest integer
such that $\liminf_q \frac{a_q}{q^k} > 0$ and call this number the s-dimension
of $R$, $\sdim(R)$.  We can then prove 

\begin{Thm} There exists a prime ideal
$\p = \p(R)$ such that $R/P$ is strongly $F$-regular and $\sdim(R) \leq
\dim(R/\p)$.  The prime $\p$ is $0$ if and only if $R$ is strongly $F$-regular.
\end{Thm} 

The proof of this theorem is contained in Theorem~\ref{prime} and
Theorem~\ref{main}.

The first part of the paper will introduce the relevant concepts of the paper,
such as the splitting dimension and the Frobenius splitting ratios. 
Section  2 of the paper contains some results about the splitting dimension, the
splitting prime ideal $\p$, and the lower and upper Frobenius splitting ratios of a ring
$R$. We also raise the question of whether or not $\sdim(R) = \p(R)$ and provide some
evidence in support of an affirmative answer. The study of these
concepts continues 
in Section 3 where we investigate
the case of reduced rings that are images of regular local rings and state our
main Theorem (Theorem~\ref{main}). The last section of the paper extends some of
the earlier results to the class of $R$-modules.

We would like to conclude our introduction by outlining some of the main facts
about $F$-purity and strong $F$-regularity. The notion of $F$-purity was
introduced in the work of Hochster and Roberts (\cite{HR}) and that of strong
$F$-regularity in a paper by Hochster and Huneke (\cite{HH-str}). Both concepts
have been studied by many other authors;  the main facts about them,
which are now considered standard, will be listed below without any proofs or
references.

Let $F:R \to R$ be the Frobenius homomorphism $F(r)=r^p$. We denote by
$F^e$ the $e$th iteration of $F$, that is $F^e(r) = r^{q}$, $F^e:R
\to R$. One can regard $R$ as an $R$-algebra via the homomorphism
$F^e$. Although as an abelian group it equals $R$, it has a different
scalar multiplication. We will denote this new algebra by $R^{(e)}$.

\begin{Def}
 $R$ is \emph{$F$-finite}  if $R^{(1)}$ is module finite over $R$, or, 
equivalently (in the case that $R$ is reduced),
 $R^{1/p}$ is module finite over $R$.  R is called 
\emph{$F$-pure} if  the Frobenius homomorphism is a pure map, i.e,  
 $F \otimes_R M$ is injective for every $R$-module $M$.
\end{Def}

If $R$ is $F$-finite, then $R^{1/q}$ is module finite over $R$, for
every $q$. Moreover, any quotient and localization of an $F$-finite
ring is $F$-finite. Any finitely generated algebra over a perfect
field is $F$-finite. An $F$-finite ring is excellent.

If $R$ is $F$-finite, then $R$ is $F$-pure if and only if the inclusion $R \hookrightarrow R^{1/p}$ splits over $R$. This is equivalent to the assertion that $R \hookrightarrow R^{1/q}$ splits, for every $q$. The reader should note that $R^{1/q} \simeq R^{(e)}$ as $R$-algebras, by sending $r^{1/q}$ to $r$.

Let $R^o$ be the complement of the minimal primes of the ring $R$.

\begin{Def}
A reduced Noetherian $F$-finite ring $R$ is \emph{strongly $F$-regular}  if for every $c \in R^o$ there exists $q$ such that the $R$-linear map $R \to R^{1/q}$ that sends $1$ to $c^{1/q}$ splits over $R$, or equivalently $Rc^{1/q} \subset R^{1/q}$ splits over $R$.
\end{Def}

The notion is related to that of $F$-purity by the fact that, in the definition above, if a splitting exists for a choice of $c \in R^o$ and $q$, then the inclusion $R \hookrightarrow R^{1/q'}$ (sending $1 \mapsto 1$) also splits for every $q'$. Also, if a splitting of $Rc^{1/q} \subset R^{1/q}$ exists, then a splitting of $Rc^{1/q'} \subset R^{1/q'}$ exists too for all $q' \geq q$.

The notion of strong $F$-regularity localizes well, and all ideals are tightly closed in strongly $F$-regular rings. Regular rings are strongly $F$-regular and  strongly $F$-regular rings are Cohen-Macaulay and normal. To show that a reduced Noetherian and $F$-finite ring is strongly $F$-regular it is necessary and sufficient to show that, for some $c \in R^o$ such that $R_c$ is strongly $F$-regular, there exists $q$ such that $Rc^{1/q} \subset R^{1/q}$ splits over $R$.

Let $E_R(k)$ be the injective hull of the residue field of $R$. Then $R$ is strongly $F$-regular if and only if $\et =0$. The ideal $\ann_R (\et)$ is called the $CS$-test ideal of $R$ and is denoted by $\cs (R)$. Note that $\cs (R) = R$ if and only if $R$ is strongly $F$-regular. 

\begin{Def}\label{approx-gor}
The ring $(R,\fm)$ is \emph{ approximately Gorenstein} if
$R$ has a sequence of $\fm$-primary irreducible ideals $\{I_t\}_t$ cofinal
with the powers of $\fm$. See \cite{Ho}.
\end{Def}
 By taking a subsequence, we may assume that
$I_t \supset I_{t+1}$.  For each $t$, let $u_t$ be an element of $R$ which
represents a socle element modulo $I_t$.  Then there is, for each $t$, a
homomorphism $R/I_t \hookrightarrow R/I_{t+1}$ such that $u_t +I_t\mapsto u_{t+1} + I_{t+1}$.  The direct limit of the system will be  the injective hull $E =E_R(R/\fm)$
and each $u_t$ will map to the socle element of $E$, which we will denote by $u$.  Hochster has shown that every excellent, reduced local ring is
approximately Gorenstein (\cite{Ho}).  In particular, every $F$-finite reduced ring is approximately Gorenstein.


\section{Notation and terminology}

Let $\ringR$ be a reduced, local, $F$-finite ring of positive characteristic $p>0$
and Krull dimension $d$. As in the introduction, let $$R^{1/q} = R^{a_q} \oplus M_q$$ 
be a
direct sum decomposition of $R^{1/q}$ such that $M_q$ has no free
direct summands. If $R$ is complete, such a decomposition is unique up
to isomorphism, but it is always the case that the values of the $a_q$ are
independent of the decomposition. As mentioned in the introduction, we investigate here the size of the
numbers $a_q$ with respect to $q$. 

For a local ring $(R,\fm , k)$,
we set $\alpha(R) =\log _p[k:k ^p]$. For $I$ an ideal of $R$, let $I^{[q]} = (i^q: i\in I)R$.
 It is easy to see that, for
an $\fm$-primary ideal $I$ of $R$, $\length_R (R^{1/q}/IR^{1/q})=
\length(R/I^{[q]}) q^{\al(R)}$. If $P$ is a prime ideal of $R$, then
$\al(P)$ will simply denote $\al(R_P)$.  It is often convenient to
understand many of the statements in this paper first in the case that $k$ is perfect, i.e., $\al(R) = 0$.  However, in order to use induction and have the concepts localize
well we must keep track of the exponent $\al(R)$.

Before we proceed, we would like to remind the reader of the following
fact (\cite{K}).

\begin{Lma}
Let $R$ be an $F$-finite Noetherian ring of characteristic $p$. Then
for any prime ideals $P \subset Q$ of $R$, $$\al(P) = \al (Q) +
\height(Q/P).$$
\end{Lma}

We would like to first define the notion of $F$-signature as it appears in~\cite{AL} and~\cite{HL}.

\begin{Def}
 Let $R$ be an $F$-finite local ring.
The \emph{$F$-signature} of $R$ is $$s(R) = \lim_{q \to \infty}
a_q/q^{d+\al(R)},$$ 

\noindent
if it exists.
\end{Def}

The following result, due to Aberbach and Leuschke, holds:

\begin{Thm}\label{AL}
Let $\ringR$ be a reduced Noetherian ring of positive characteristic $p$. Then
$\liminf_{q \to \infty} a_q/q^{d+\al(R)} >0$ if and only if $\limsup_{q \to \infty} a_q/q^{d+\al(R)} >0$ if and only if $R$ is strongly $F$-regular.
\end{Thm}

Moreover, Huneke and Leuschke proved that if $R$ is local, reduced and Gorenstein, then the $F$-signature exists. Yao has recently extended this result to rings that are Gorenstein on their punctured spectrum.

We are now able to define the central concepts to be investigated in this paper.

\begin{Def} 
 Let $R$ be an $F$-finite local ring.
The \emph{$s$-dimension}  of R, $\sdim(R)$, is the largest integer $k$
such that $$\liminf _{q \to \infty} \frac{a_q}{q^{k+\al(R)}}$$ is not
zero. If $R$ is not $F$-pure then $\sdim(R) = -\infty.$ The \emph{lower Frobenius splitting ratio} of 
$R$, $\frm(R)$,  equals
the value of the limit introduced above, that is  $\frm(R) = \liminf
_{q \to \infty} \frac{a_q}{q^{\sdim(R)+\alpha(R)}}$. Similarly,  the
\emph{upper Frobenius splitting ratio} of $R$, $\frp(R)$, equals $\limsup
_{q \to \infty} \frac{a_q}{q^{\sdim(R)+\alpha(R)}}$. Whenever $\frm
(R) = \frp (R)$, we call this number the \emph{Frobenius splitting ratio} of
$R$, $r_F(R)$.

 Whenever $\sdim(R) = \dim (R)$ (equivalently, $R$ is strongly $F$-regular), we recover the notion of $F$-signature of $R$, denoted by $s(R)$, and, we will use this terminology whenever we are in such a special case. In particular,  whenever $\sdim(R) =\dim (R)$, we will be speaking of the \emph{lower and upper signature} of $R$,  denoted accordingly by $\spl(R)$ and $\smi(R)$, instead of its lower and, respectively, upper Frobenius splitting ratio.
\end{Def}

We show below that $R$ is $F$-pure if and only if $\sdim(R) \geq 0$.

\begin{Rem}
 In \cite{WY} the notion of minimal Hilbert-Kunz multiplicity  is defined, denoted $m_{HK}(R)$.  Theorem~\ref{AL} shows that this number is non-zero
if and only if $R$ is strongly $F$-regular.  In our terminology, $m_{HK}(R) = 0$ if $\sdim(R) < \dim(R)$ and $m_{HK}(R) = s^-(R)$ otherwise.
\end{Rem} 

\begin{excise}
{\begin{Rem}{\rm Our definition of $F$-signature differs from the one
existent in the literature where the $F$-signature of $R$ equals
$\liminf _{q \to \infty} \frac{a_q}{q^{d+d}}$. However, it is known
that this number is nonzero if and only if and $R$ is strongly
$F$-regular (Aberbach-Leuschke). Extending the definition such that
the $F$-signature is not trivial for other  classes of rings than the
strongly $F$-regular ones will prove to be both useful and natural
later in the paper; no information is lost this way, since the
exponent of $q$ that ``gives'' the $F$-signature is described by the
$s$-dimension of $R$.}
\end{Rem}
}
\end{excise}

The following result is classical and very helpful in
understanding the number of free direct summands of $R^{1/q}$. We will quote it here in the form stated in~\cite{HH-adv}.

\begin{Thm}
\label{hh}
Let $\ringR$ be a local ring and $F$ be a free finitely generated
$R$-module. Let $M$ be a finitely generated $R$-module and $\phi: F
\to M$ an  $R$-linear map. Then $\phi$ splits over $R$ if and
only if $\phi \tensor E_R(k)$ is an injective map.
\end{Thm}

If $M$ is a finite $R$-module that has no nonzero free direct
summands, then for every $m \in M$, the map $R \to M$ that sends $1$
to $m$ does not split. By tensoring with $E_R(k)$ and using the above
Theorem, we conclude that the map $E_R(k) \to E_R(k) \tensor M$ that
sends $e$ to $e \tensor m$ is not injective. Denote by $u$ the socle
generator for $E_R(k)$. The following holds:

\begin{Prop} 
\label{socle}

Let $\ringR$ be a local ring, let $M$ be a finite $R$-module, and let $u$ be the socle generator 
of $E_R(k)$. Then $M$
has no nonzero free direct summands if and only if for every $m \in M$,
$u \tensor m=0$
in $E_R(k) \tensor M$.

\end{Prop}

\begin{proof}
For each $m \in M$, define $f_m : R \to M$ by $f_m(r)=rm$. Then $f_m$does not split over $R$ if and only if $f_m
\tensor E_R(k)$ is not injective if and only if ${\rm ker}(f_m \tensor
E_R(k))$ is nonzero (Theorem~\ref{hh}). This holds if and only if
${\rm ker}(f_m \tensor E_R(k))$ intersects the socle of $E_R(k)$
non trivially. But, $\Soc(E_R(k)) = Ru$. So, $f_m \tensor E_R(k)$ is
not injective if and only if $u \tensor m =0$
\end{proof}

The following Corollary will prove to be of much use in what follows. Yao has
also noted this result in his recent work~\cite{Y}.

\begin{Cor}
\label{aq} 
Let $\ringR$ be an $F$-finite reduced local ring of characteristic $p$. Let $R^{1/q} =
R^{a_q} \oplus M_q$ be a direct sum decomposition of $R^{1/q}$ over
$R$, where $M_q$ has no $R$-free direct summands. Define $$J_q: = \{r \in R^{1/q} : r \tensor u =0 \ \ {\rm in} \ \
R^{1/q} \tensor E_R(k) \}.$$ Then, $a_q = \length_R (R^{1/q}/J_q).$

\end{Cor}
\begin{proof}
$R^{1/q} = R^{a_q} \oplus M_q$; so, by tensoring with $Ru \subset
E_R(k)$, we get $$R^{1/q} \tensor R \cdot u = (R^{a_q} \tensor R \cdot
u) \oplus (M_q \tensor R \cdot u).$$ Map this further to $R^{1/q}
\tensor E_R(k)$ and apply Proposition~\ref{socle} to $M_q$. We have
that $\im(R^{1/q} \tensor R \cdot u \to R^{1/q} \tensor E_R(k) ) = \im
(g:R^{a_q} \tensor Ru \to R^{1/q} \tensor E_R(k))$, where $g(a \tensor
ru) = a \tensor ru$. But $g$ factors through $R^{a_q} \tensor E_R(k)
\to R^{1/q} \tensor E_R(k)$ which is an injective map (again by
Theorem~\ref{hh}), so $g$ is also injective.  But $\length(Ru) =1$,
and hence $\length(R^{a_q} \tensor Ru)=a_q.$ In conclusion,
$$\length(\im(R^{1/q} \tensor R \cdot u \to R^{1/q} \tensor E_R(k) ))
= a_q.$$ There is a natural $R$-linear map $R^{1/q} \to R^{1/q}
\tensor E_R(k) $ that sends $r \in R^{1/q}$ to $r \tensor u$.  The
kernel of this map equals $J_q$ and the theorem is proved.
\end{proof}

\section{Splittings, $s$-dimension and  strong $F$-regularity}

Let $F^e : R \to R$ be the $e$th iteration of the Frobenius map. We
can regard $R$ as a new (right) $R$-algebra, under this map, which we
denote by $R^{(e)}$.
Corollary~\ref{aq} leads us to define the following sequence of ideals.

\begin{Def} 
\label{Aq}  Let $(R,\fm,k)$ be an $F$-finite local ring.
  We define $$A_e = A_e (R) = \{ r \in R: r \tensor u =0
\in R^{(e)} \tensor_R E_R(k) \}.$$  
\end{Def}
If $R$ is not reduced then the map $R \to R^{(e)}$ cannot be injective
for $e >0$, and hence by Theorem~\ref{hh}, $A_e = R$ for all $e \ge 1$.
So assume that $R$ is reduced (hence approximately Gorenstein).
Then clearly, $A_e = (J_q)^{[q]}$ (note that 
$R^{1/q} \simeq R^{(e)}$ by $r^{1/q} \longmapsto r$). In fact, $A_e$ is an ideal
both of $R$ and $R^{(e)}$. 
Using the notation in Definition~\ref{approx-gor}, we see that 
$A_e = \cup_t (I_t\brq:u_t^q)$.  We claim that the sequence $\{A_e\}_e$
is non-increasing. If $R$ is not $F$-pure then $A_e = R$ for
all $e$.  Otherwise, let $r \in A_{e}$.  For $t \gg 0$, $r u_t^{q} \in I_t^{[q]}$ where $q = p^e$, hence
$r^p u_t^q \in I_t\brq$.  Taking $p$th roots and applying the splitting
yields, $r u_t^{q/p} \in I_t^{[q/p]}$, showing that $A_e \subseteq A_{e-1}$.

\begin{Def}  Let $(R,\fm,k)$ be an $F$-finite reduced local ring.
 Let   
$$
\p(R) =  \{ r \in R: r \tensor u =0 \in
R^{(e)} \tensor_R \E \text{  for  all\ } e \gg 0 \} = \cap_e A_e.
$$  
If the ring is understood we will sometimes refer to this ideal
simply as $\p$.
\end{Def}
 As
we will show in Theorem~\ref{prime} below, the ideal $\p(R)$ is prime.

Since $u$ is the socle generator for $\E$, we have $r^q \otimes u = 1
\otimes ru$, which is $0$ if $r \in \fm$, so  $m^{[q]} \subset A_e$ for every $q$. On
the other hand, if $u \notin \et$, then $\p =0$ (see~\cite{A}, Proposition 2.4). If $R$ is a domain, then $\p$ is
nonzero only if $R$ is not strongly $F$-regular, since $ \et = 0$ in strongly $F$-regular rings (see
also Corollary~\ref{ineq}). 
The reader is invited to note that there is a similarity between the ideal
$\p(R)$ introduced here and the concept of $F$-stable primes of $R$
introduced in~\cite{E}. In fact, when $R$ is Gorenstein, then $R$
admits a unique $F$-stable prime, which, indeed, is $\p$. This
similarity is underlined by our next Theorem that asserts that $\p(R)$
is in fact a prime ideal of $R$.  We will refer to $\p(R)$ as the
\emph{splitting prime} for $R$.

Before we state the result, we need to introduce some notation. For
every $c \in R$, one can define an $R$-linear map $\phi _{c,e} : R \to
R^{1/q}, \phi _{c,e} (1) = c^{1/q}$, where $q=p^e$.

\begin{Thm}
\label{prime}
Let $\ringR$ be an $F$-finite reduced
 local ring of characteristic $p$. Then, $A_e = \{c
\in R: \phi_{c,q}$ does not split over $R \}$ and $\p (R)$ is a prime
ideal or the unit ideal.
\end{Thm}

\begin{proof}
By Theorem~\ref{hh} we know that $\phi_{c,e}$ splits if and only
$\phi_{c,e} \tensor \E$ is injective.

Take $c \in A_e$. Then $c^{1/q}$ belongs to $J_q$. Hence $c^{1/q}
\tensor u =0$, and therefore $\phi _{c,e}$ does not split. Now, assume
that $c \in R$ is an element such that the map $\phi_{c,q}$ does not
split. Then again $u \in {\rm Ker}(\phi_{c,e} \tensor \E)$ which means
that $c^{1/q} \in J_q$. So, $c \in A_e$. This shows the first
assertion of the Theorem.

For the second assertion assume that $\p(R)$ is not the unit ideal and let 
us write $$\p = \{ c \in R: \phi_{c,e}
\text{ does  not split for  any}\ e \gg 0 \}.$$ If
$\p$ is not prime, then there exist $a$ and $b$ not in $\p$ such that
$ab \in \p$. Since $a$ and $b$ are not in $\p$ we can find $e_1$ and
$e_2$ large enough such that the $R$-linear maps $\phi _{a,e_1}$ and
$\phi_{b,e_2}$ split. Let $\theta_1:R^{1/q_1} \to R$ an $R$-linear
that maps $a^{1/q_1}$ to $1$. This map exists because $\phi_{a,e_1}$
splits over $R$. Similarly, we can find $\theta_2:R^{1/q_1 q_2} \to
R^{1/q_1}$, such that $\theta_2 (b^{1/q_1 q_2}) =1$ and $\theta_2$ is
$R^{1/q_1}$-linear (we can do this by taking the $q_1$-root of a
splitting for $\phi_{b,e_2}$).

Let $\theta = \theta_1 \theta_2$. It is easy to check that $\theta$ is
$R$-linear and $\theta ((a^{q_2}b)^{1/q_1q_2})=1$. So, for $e_0 =e_1
+e_2$, $\phi_{a^{q_2}b, e_0}$ splits.

We know that $ab \in \p$, and this implies that $a^{q_2}b \in \p$. So,
if we fix $e_2$, this means that $\phi_{a^{q_2}b,e}$ does not split
for any $e \gg 0$ large enough. This is a contradiction, since we have
just shown the existence of a splitting for values $q$ of the form
$q=q_0 =p^{e_0}=p^{e_1+e_2}$, where $q_1=p^{e_1}$ can be arbitrarily large.

In conclusion, $\p $ is a prime ideal.
\end{proof}

\begin{excise}
{Ultimately, we will show that $\sdim(R) = \dim(R/\p)$ (see Theorem~\ref{main}).
For the moment we can show one inequality.}
\end{excise}

\begin{Cor}
\label{ineq}
Let $\ringR$ be an $F$-finite reduced  local  ring of characteristic $p$.  Then, $R$ is $F$-pure if and only if $\p(R)$ is a proper ideal if and only if
$\sdim(R) \ge 0$.  In this case,
$$
\sdim(R) \leq \Dim(R/\p).$$  The ring $R$ is strongly $F$-regular if
and only if $\p  = 0$. Also, $\Dim(R/\p)=0$ if and only if $\p = \fm$
if and only if $\sdim(R) = 0$, in which case $\fr(R) = 1$.
\end{Cor}

\begin{proof}
$R$ is $F$-pure if and only if some $a_q > 0$ if and only if $A_e$ is
proper, and since the sequence $\{A_e\}_e$ is non-increasing this occurs
if and only if $\p(R) \ne R$.  By Theorem~\ref{aq}, $a_q = \length (R^{1/q}/J_q)$
and hence $a_q = \length(R/A_e) \cdot q^{\al(R)}$ with $\alpha(R) =\log _p[k:k ^p]$ .   From this we can see that $R$ is $F$-pure if and only if 
$\sdim(R) \ge 0$.

  We next show that $\sdim(R) \leq \Dim(R/\p).$ Since $\fm^{[q]} \subset A_e$, for every $q = p^e$, we have that $\fm ^{[q]}
+ \p \subset A_e$.    Then
$a_q/q^{k+\al(R)} = \length(R/A_q)/q^k \leq
\length(R/\fm^{[q]}+\p)/q^k$.  If $k > \dim(R/\p)$, then by the theory
of Hilbert-Kunz functions, $\lim_{ q \to \infty}
\length(R/\fm^{[q]}+\p)/q^k =0$.  So, $\sdim(R) \leq \Dim(R/\p)$.

The ideal $\p = 0$ if and only if for all $c \ne 0$ there exists a $q$
such that $\phi_{c,q}$ splits; this is precisely the definition of a
strongly $F$-regular ring.

Since $\p$ is prime,   $\Dim(R/\p) =0$ if and only if
$\p=\fm$. Clearly $\p \ne R$ if and only if $R$ is $F$-pure if and only if $A_e \neq R$ for all $e$. So,
$\p = \fm$ if and only if $R$ is $F$-pure and $A_e = \fm$ for all $e$.
In conclusion, $a_q/q^{\al(R)}= \length(R/A_e) = \length(R/\fm) =1$. 
Hence, in
particular, $\fr(R)=1$.
\end{proof}

\begin{Rem}  We note the following facts:

\begin{itemize}
\item[1)]  The proof of the above Corollary shows that, in fact, when $R$ is
   $F$-pure and $\p=\fm$ (an important such case is that of an
   isolated non-strongly $F$-regular point), then $a_q = q^{\al(R)}$.
\item [2)]Assume that $R$ is strongly $F$-regular. A result of Aberbach
(\cite{A}) shows that there exists $e_o$ such that $A_e \subset
\fm^{[q/q_o]}$ for every $e \geq e_o$. This shows that
$a_q/q^{d+\al(R)} \geq \length(R/\fm^{[q/q_o]})/q^d$; since
$e_{HK}(\fm) >0$, we conclude that $\sdim(R)=d$ and $\fr(R) \geq
e_{HK}(\fm)/q_o^d$. This result was first established in~\cite{AL} and
our remark is along the lines of the argument given there.
\item [3)] Also, let us take $c \in \cs(R) =\ann_R(\et)$, a $CS$-test element. The ring $R_c$
is strongly $F$-regular. If $c \notin A_e$, for some $e$, then $R$ is strongly $F$-regular (as explained in the introduction). So, if $R$ is not strongly $F$-regular, then $\cs(R) \subseteq \p$.
\end{itemize}
\end{Rem}

The splitting prime of $R$ does not localize well in general.  The most obvious
example is when $\p(R) = \mathfrak m$.  Then for all $Q \in \Spec(R)$, $R_Q$ is
$F$-pure, so $\p(R_Q) \ne R_Q = \p(R)R_Q$.  However, when $Q \supseteq \p(R)$
we do get a nice localization result.

\begin{Prop}
Let $\ringR$ be a local ring of characteristic $p$. Then for every $Q
\in \Spec(R)$ such that $\p \subset Q$ we have
$$ \p(R_Q) = \p(R) R_{Q}.$$
\end{Prop}

\begin{proof}
We first note that $(R^{1/q})_Q = (R_Q)^{1/q} = R^{1/q} \tensor
R_Q$. Moreover, if $\phi _{c,e}$ splits, then $\phi_{c,e} \tensor R_Q$
splits, too.

Take $c/s \in \p (R_Q)$. Then $c/1$ belongs to $\p (R_Q)$. We would
like to show that $c \in \p(R)$. Assume the contrary; then,
$\phi_{c,e}:R \to R^{1/q}$ does split for $e \gg 0$. By tensoring with
$R_Q$ we get a splitting over $R_Q$, which is a contradiction.

Now, take $c \in \p(R)$. We want to show that $c/1$ belongs to $\p
(R_Q)$. If not, for large $q$, the map $R_Q \to R^{1/q}_Q$ that sends
$1$ to $c^{1/q}/1$ splits. So, for each such $q=p^e$, there is an element
$d_q \notin Q$ and a map $R^{1/q} \to R$ that sends $c^{1/q}$ to
$d_q$. But, $d_q$ is not in $\p(R)$, since $\p(R) \subset Q$. So, for
$q' \gg 0$, there is an $R$-linear map $R^{1/q'} \to R$ that sends
$d_q^{1/q'}$ to $1$. Combining these last two assertions, we conclude
that $c$ is not in   $\p(R)$ and this is a contradiction.

\end{proof}

\begin{Prop}
Let $R$ be an $F$-pure local ring. Then $Q \subset \p$, for every $Q
\in {\rm Min}(R)$.

\end{Prop}

\begin{proof}

Let $c \in R$. For every $q$, the map $\phi_{c,q}$ is not injective if
$c$ is a zero-divisor on $R$. Since $R$ is reduced this means that $c$
belongs to some minimal prime $Q$ of $R$. So, if $c \in Q$, with $Q$ a
minimal prime, then $\phi_{c,q}$ cannot split over $R$ for any
$q$. The assertion of the Proposition follows.

\end{proof}

We would like to show that the splitting prime of the completion is
extended.  In order to do this we need the next Lemma.
\begin{Lma}\label{flatmap}
Let $(R,\fm) \to (S,\fn)$ be a faithfully flat map of $F$-finite local reduced
rings such that the closed fiber $S/\fm S$ is regular of dimension $n$.  Then
$A_e(S) = A_e(R)S + \fn \brq$. If $\sdim R =0$, then $\sdim S = n$, $\p(S) = \fm S$, and $\frm(S)=\frp(S) = 1$.
\end{Lma}

\begin{proof}
Both rings are approximately Gorenstein, and the closed fiber $S/\fm S$ is regular of dimension $n$. Let $z_1,...,z_n$ be a sequence of elements in $S$ that define a regular sequence on $S/\fm S$.  Then for any $q=p^e$, $z_1^q,...,z_n^q$ form a regular sequence on $S/\fm S$ and, in fact, $R \to S/(z_1^q,...,z_n^q)S$ if a faithfully flat map.

If $\{I_t\}_t$ is a sequence of irreducible $\fm$-primary
ideals of $R$, cofinal with the powers of $\fm$, then $\{J_t \}_t=\{ I_t S +(z_1^t,...,z_n^t) \}_t$ is
a sequence of irreducible $\fn$-primary
ideals of $S$, cofinal with the powers of $\fn$.  Let $u_t$ be elements
of $R$ which are socle elements modulo $I_t$.  Then the images of the $u_t \cdot (z_1\cdots z_n)^{t-1}$
in $S$ are socle elements modulo $J_t$.  Moreover, $R/I_t \brq  \to S/I_t \brq S$ is flat with regular closed fiber of dimension $n$. Since $R/I_t \brq$ has depth $0$,  we see that $S/I_t \brq  S$ has depth $n$. It follows that $S/ I_t \brq S$ is Cohen-Macaulay. In particular, this shows that
$z_1,...,z_n$ form a regular sequence on $S/I_t \brq S$, because their images in $S/I_t \brq S$ form a system of parameters. 

Let us compute $A_e(S) = \cup_t ((I_t \brq S +(z_1^{tq}, \cdots, z_n^{tq})) : u_t^q \cdot (z_1 \cdots z_n)^{(t-1)q})$. If $a \in S$ such that

$$ a \cdot u_t^q \cdot (z_1 \cdots z_n)^{(t-1)q}   \in I_t \brq S +(z_1^{tq}, \cdots, z_n^{tq})$$
\noindent
then, by using that $z_1,...,z_n$ form a regular sequence on $S/I_t \brq S$, we have that  

$$a u_t^q \in I_t \brq S +(z_1^q,...,z_n^q).$$
\noindent
This is equivalent to $a \in (I_t \brq S + (z_1^q,...,z_n^q)):_S u_t^q)$. 

By  flatness of $S/(z_1^q,\cdots, z_n^q)$ over $R$,
$(I_t \brq S +(z_1^q,...,z_n^q)):_S u_t^q = (I_t \brq :_R u_t^q)S +
(z_1^q,...,z_n^q)$. In conclusion, $A_e(S) = \cup_t ((I_t \brq :_R
u_t^q)S + (z_1^q,...,z_n^q))$ and hence $A_e(S) = A_e(R)S
+(z_1^q,...,z_n^q)$.  So, if $\sdim(R)=0$, then $A_e(S)= \fm \brq S
+(z_1^q,...,z_n^q)$. Hence, $a_q = \length_S (S/ A_e(S))= \length _S (S/(\fm S +(z_1^q,...,z_n^q))$ which shows that $\sdim(S)=n $.
\end{proof}

\begin{Prop}
\label{completion}
Let $(R, \fm)$ and $(S, \fn)$ be  local, reduced, $F$-finite rings of characteristic $p$ and let $f : (R, \fm)  \to (S, \fn)$ be a flat local map with regular fibers.  Then $A_e (S) = A_e(R) S + \fn \brq$ for every $q = p^e$, and $\p (S) = \p(R) S$. In particular, if $R$ is local, reduced and $F$-finite, then 
$A_e (R) \widehat{R} = A_e(\widehat{R})$ and $\p(R) \widehat{R} = \p(\widehat{R})$.
\end{Prop}

\begin{proof} The assertion that $A_e(S) = A_e(R) + \fn\brq$ has been
shown in Lemma~\ref{flatmap}.

To see that $\p(S) = \p(R) S$ we first observe that
 $\p(S) = \cap_e A_e(S)$ is a prime ideal of $S$
lying over $\p(R)$.  

We first show that any prime of $S$ that is minimal over $\p(R)S$ is in fact contained in $\p(S)$.   Let
$Q  \ne \p(S)$ be a prime ideal of $S$ that is minimal
over $\p(R)S$. Since $R \to S$ is flat then $Q$ lies over $\p$.  Thus the map $R_{\p(R)} \to S_Q$ is faithfully
flat and $\p(R) S_Q = Q S_Q$.  By Lemma~\ref{flatmap}, $\sdim (S_Q) =0$.  
In particular, $A_e(S_Q) = Q S_Q$ for all $e$.
 However, if we
take $c \in Q - \p(S)$, then for $q \gg 0$,  $c^{1/q} {S} \subseteq {S}^{1/q}$ splits over $S$.  This remains true after
localizing, showing that $A_e(S_Q) \ne Q S_Q$, a
contradiction.

In conclusion, $Q \subseteq \p(S)$.

Note that $\p(R)$ is prime in $R$, and since $R\to S$ has regular fibers, then
$\p(R)S$ is a radical ideal in $S$.

Now, we would like to argue that $\p(S)$ itself is a minimal prime of  $\p(R)S$. The homomorphism $R_{\p(R)} \to S_{\p(S)}$ is flat because $\p(S)$ lies over $\p(R)$. Its closed fiber is regular and $\sdim(R_{\p(R)}) =0$. If the dimension of its closed fiber is positive, then, by Lemma~\ref{flatmap} below, $\sdim(S_{\p(S)})>0$. However, this is impossible. In conclusion, the closed fiber of $R_{\p(R)} \to S_{\p(S)}$ is 0-dimensional, so  $\p(S)$ is minimal over $\p(R) S$.

Now, since $\p(R)S$ is radical and its only minimal prime is $\p(S)$, we obtain the equality
$\p(R)S = \p(S)$.

The last assertion about the completion follows from the fact that, under our hypotheses, $R$ is excellent, so all the fibers of the completion homomorphism are smooth, hence regular.
\end{proof}

 Proposition~\ref{completion} enables us to reduce the study of the $s$-dimension
and the Frobenius splitting ratio of reduced excellent rings to the case of complete local
rings. These rings are images of regular local rings.  Following in the
paths of Fedder,  Glassbrenner, and Cowden-Vassilev we can analyze the
finer structure of $F$-pure rings in this way.

\section{Images of regular rings}

In this section, we continue our analysis in the case of local rings
that are images of regular local rings. The more general case of
reduced F-finite rings can be
reduced to this class of rings by Proposition~\ref{completion}. 

Let us write $R=S/I$ where $S$ is regular and local and $\pi : S \to
R$ the natural projection. Assume that $R$ is reduced, i.e., $I$ is
a radical ideal. There are two results that describe under  what conditions
such rings are $F$-pure or strongly $F$-regular. The first one is due to
Fedder (Theorem 1.12 in~\cite{F}).  It does not require $R$ to be $F$-finite.

\begin{Thm}
Let $(S, \fm, k)$ be a regular local ring of characteristic $p$ and
let $R=S/I$. Then $R$ is $F$-pure if and only if $(I^{[p]}:I)
\not\subset \fm ^{[p]}$.
\end{Thm}

The idea behind this criterion for $F$-purity has been used by
Glassbrenner to give a similar criterion for strong $F$-regularity
(Theorem 2.3 in~\cite{G}). We will state it here in the form we need
it later.

\begin{Thm}
\label{G}
 Let $(S, \fm, k)$ be a local regular ring and let $R=S/I$ and $c \in
S$. Then, for every $q=p^e$, the map $\phi_{\pi(c),e}$ splits if and only
if $ c \notin \fm^{[q]} : _S (I^{[q]}:I)$.  Moreover,  $R$ is strongly
$F$-regular if and only if $I = \cap_q \bigl(\fm^{[q]} : _S (I^{[q]}:I)\bigr)$.
\end{Thm}

Using the notation introduced previously, we see that $A_e =
\bigl(\fm^{[q]} : _S (I^{[q]}:I)\bigr)/I.$ So, $a_q = q^{\al(R)}
\length_S(S/(\fm^{[q]} : _S (I^{[q]}:I)))$. Note that $S/\fm^{[q]}$ is
zero dimensional and Gorenstein.   The quotient $S/\bigl( \fm^{[q]} : _S
(I^{[q]}:I) \bigr) $ is an $S/\fm^{[q]}$-module and its Matlis dual is  $
\frac{(I^{[q]}:I) +\fm ^{[q]}}{\fm^{[q]}}$. Hence, $a_q = q^{\al(R)}
\length_S(\frac{(I^{[q]}:I) +\fm ^{[q]}}{\fm^{[q]}})$.

Although Theorem~\ref{G} was stated in~\cite{G} under the assumption that
$S$ is local, regular and $F$-finite, one can note that the $F$-finiteness
hypothesis can be removed along the lines of the argument used by
Fedder in proving his criterion for $F$-purity (where we interpret ``strongly $F$-regular'' to mean that $0^*_E = 0$). The point is that one can make a flat
base change, by enlarging the residue field to its perfect closure,
to get to the case 
where $S$ is $F$-finite. Strong $F$-regularity
commutes with this base change (see, for example, Theorem
3.6~\cite{A}).

\begin{Rem}[for proofs, see \cite{F}, 1.4 and 1.5]
\label{fedder}
 For every $q=p^e$, $\Hom{S}{S^{1/q}}{S}$ has an
$S^{1/q}$-~structure given by $s^{1/q} \cdot \phi (t^{1/q}) := \phi
((st)^{1/q})$, for every $s,t \in S$ and $\phi \in \Hom
{S}{S^{1/q}}{S}$. Moreover, $\Hom{S}{S^{1/q}}{S} \simeq S^{1/q} $. Let
$T$ be a generator for $\Hom{S}{S^{1/q}}{S}$ and $s ^{1/q} \in
S$. Then $s^{1/q} T$ defines an element in $\Hom{R}{R^{1/q}}{R}$ if and only if $s^{1/q}
\in (IS^{1/q}:I^{1/q})$ if and only if $s \in (I^{[q]}:I)$.
\end{Rem}

The following result is similar to Theorem 3.1 in \cite{C}, in which Cowden-Vassilev showed that if $R$ is $F$-pure then so is $R/\tau(R)$, where
$\tau(R)$ is the test ideal of $R$.  We give
two proofs: one is similar to the proof of Theorem 3.1 given in \cite{C}, while the
other one uses some ideas originating in Fedder's work, \cite{F}.

\begin{Prop}
\label{sdim}

Let $(S, \fm, k)$ be an $F$-finite local regular ring and let $R=S/I$ be reduced
and $F$-pure. Set $\q$ to be the full preimage of $\p(R)$ in $S$. Then
for all $q$, $(I^{[q]}:I)  \subset (\q^{[q]}:\q)$.  In particular, if
$R$ is $F$-finite then 
$\sdim(R/\p) \geq \sdim(R).$
\end{Prop}

\begin{proof}
{\it Proof 1}.  $R$ is excellent and reduced, so it is approximately 
Gorenstein.  Let $\{I_t\}$ be a collection of irreducible
$\fm$-primary ideals in $R$ cofinal with the powers of $\fm$, and
denote by $\{J_t\}$ their full preimages in $S$.  
Let $u_t$ be an element of $S$
mapping to the socle in $S/J_t$.  Then we may describe $\q$ as $\q =
\cap_{q}( \cup_t(J_t^{[q]}+I):_S u_t^q)$.  Let $w \in (I^{[q]}:I)$ and
$v \in \q$.  We want to show that $vw \in \q^{[q]} = \bigl(
\cap_{q'}\cup_t(J_t^{[q']}+I):u_t^{q'}\bigr)^{[q]}$, which by flatness
of the Frobenius endomorphism over $S$ is
$\cap_{q'}\cup_t\bigl((J_t^{[qq']} + I^{[q]}):u_t^{qq'}\bigr)$.  For
all $q'$ there is a $t$ such that $vu_t^{qq'} \in J_t^{[qq']} + I$,
hence $vw u_t^{qq'} \in J_t^{[qq']} + I^{[q]}$. This shows that $vw \in
\q^{[q]}$, as desired.

{\it Proof 2}. Let $T$ be as in Remark~\ref{fedder}. In the light of this
Remark we need to show that if $\psi: = s^{1/q}T \in
\Hom{R}{R^{1/q}}{R}$, then $(s \q) ^{1/q} \in \q S^{1/q}$. By the same Remark, this is
equivalent to the assertion that $\psi$ induces an $R$-linear map
$R^{1/q}/\q^{1/q} \to R/\q$.

We can define an $R^{1/q'}$-linear map by $\psi_{q'}: R^{1/qq'} \to
R^{1/q'}$, $\psi_{q'} (a) = \psi (a^{q'})^{1/q'}$. This map is, in
particular, $R$-linear.

Take $c \in \q$; then, $c^{1/qq'} \otimes u = 0$ in $R^{1/qq'} \otimes
\E$ where $u$ is the socle generator in $\E$.

Clearly, $ \psi_{q'} (c^{1/qq'}) \otimes u = 0$. That is, $\psi
(c^{1/q}) ^{1/q'} \otimes u = 0$. So, $s^{1/q}T (c^{1/q}) = \psi
(c^{1/q})  \in \q$. Hence,  $s^{1/q} T$ takes $\q ^{1/q}$ into
$\q$. So, $s^{1/q}T$ defines an element in $\Hom
{R}{R^{1/q}/\q^{1/q}}{R/\q}$, and hence $s \in (\q^{[q]} : \q)$.

The last statement now follows, after noting that $\al(R) = \al(R/\p)$,
since 
$$
a_q(R) = q^{\al(R)} \length_S\left(\frac{I\brq:I + \fm\brq}{\fm\brq}\right) \le q^{\al(R/\p)} \length_S\left(\frac{\q\brq:\q + \fm\brq}{\fm\brq}\right) = a_q(R/\p).
$$

\end{proof}
\begin{Rem}
 This result improves Theorem~\ref{ineq}, since the $s$-dimension
of a ring is bounded above by its dimension.

\end{Rem}

\begin{Cor}
Let $R$ be a characteristic $p$ local $F$-pure  ring that is a
homomorphic image of a regular local ring $(S,\mathfrak{m},k)$, that is $R=S/I$. Then $\sdim(R_{\p})=0$.
\end{Cor}

\begin{proof}
Denote by $k(\p)$ the residue field of the localization of $R$ at
$\p$.  Apply Theorem~\ref{sdim} to $R_{\p}$. We get that
$0=\sdim(k(\p)) \geq \sdim(R_{\p})$ and we are done.
\end{proof}

\begin{Thm}
\label{str}
Let $R=S/I$ be a characteristic $p$ local $F$-pure, $F$-finite ring that is a
homomorphic image of a regular local ring $S$. Then $R/\p$ is strongly
$F$-regular.

\end{Thm}

\begin{proof}
If $\sdim R = 0$, i.e., $\p = \fm_R$ this is clear.  Thus we may assume
that $\sdim R > 0$.  
If $R/\p$ is not strongly $F$-regular we may localize a prime minimal among
the set $\{P | (R/\p)_P \text{ is not strongly $F$-regular}\}$.  After
relabeling this ring $R$, we note that the $s$-dimension is still positive.
Thus we can assume that $(R/\p)_Q$ is
strongly $F$-regular for all $Q$, different from $\fm$, containing
$\p$. We know that $R/\p$ is $F$-pure by Proposition~\ref{sdim}.

If $R/\p$ is not strongly $F$-regular, then $\p(R/\p)$ equals
$\fm_R$ since $\fm_R = \widetilde\tau(R) \subseteq \p(R)$. But then $\sdim (R/\p) =0$, a contradiction.
\end{proof}

Now we are in a position to describe the $s$-dimension of a local ring
of positive characteristic.  We also give some information on the
splitting ratios.  We expect that the upper and lower splitting ratios
always coincide.  At present, we can only prove this if $\sdim(R) \le 1$,
however, we can give an upper bound on the ratio  $r^+_F(R)/r^-_R(F)$, if $\sdim(R) = \Dim(R/\p)$.

\begin{Thm}
\label{main}
Let $R$ be an $F$-pure, $F$-finite local ring of positive characteristic. Then

i)  $R/\p$ is strongly $F$-regular, $\sdim(R) \leq \Dim
(R/\p)$. Moreover, if $\sdim(R) = \Dim(R/\p)$, then  
$\frm(R) \leq \spl (R/\p) \leq 1$, $\frp(R) \leq \smi (R/\p) \leq 1$.

ii) If $\sdim(R) = \Dim(R/\p)$, then $\frp(R) \leq \frm(R) \cdot e_{HK}(R/\p)$. This shows that, in this case, $\frp (R) > 0$ if
and only if $\frm(R) >0$.

iii) If $\Dim(R/\p) =1$, then $\sdim(R) =1 $ and $R/\p$ is a DVR.

iv) If $\Dim (R/\p) \leq 1$ then $\frp (R) = \frm (R)$.

v) $\Dim(R/\p)=\Dim(R)$ if and only if $R$ is strongly $F$-regular. In this case, $\spl (R) \leq \smi (R) e_{HK}(R)$.

vi) $\depth R \ge \sdim R$.

\end{Thm}

\begin{proof}

One can pass to the completion of $R$ by Proposition~\ref{completion}, and hence assume that $R=S/I$
with $S$ regular local. 

For part i), let us note that since $R$ is a homomorphic image of a regular ring, $R/\p$ is
strongly $F$-regular by Theorem~\ref{str}.  

\begin{excise}{First, let us remark that if $z \notin A_{e'}$, then $z^q \notin
A_{e+e'}$. Indeed, if $z^q \in A_{e+e'}$, then $z \tensor u = 0$
in $R^{(e+e')} \tensor E$. But $R$ is $F$-pure, and, so, we can
conclude that $z \tensor u=0$ in $R^{(e')} \tensor E$. This is
impossible, because $z \notin A_{e'}$.

Take $y_1,...,y_k$ parameters in $R/\p$. Let $T = L[[y_1,...,y_k]]
\subseteq R/\p$ and write $\fn$ for the maximal ideal of $T$.  Then
$\length_R(R/A_e) \geq \length_T (T/A_e \cap T)$. Let $y = y_1 \cdots
y_k$. We can find $q_o$, such that $y \notin A_{e_o}$. We claim that
$A_{e+e_o} \cap T \subset \fn^{[q]}$. If not, then the socle generator
of $T$ modulo $\fn^{[q]}$ belongs to $A_{e +e_o} \cap T$; that is,
$y^{q -1} \in A_{e+e_o}$. But, if this is true, then $y^q \in A_{e+e_o}$
which is certainly false, since $y \notin A_{e_o}$.

So, $\length_T (T/A_{e+e_o} \cap T) \geq \length_T (T/\fn^{[q]})= q^{k +
\al(R)}$. This implies that $\sdim(R) \geq k$ (and  $\frm(R) \geq
1/{q_o}^k$), and we are done.}
\end{excise}
The inequalities stated in i) follow easily from the proof of Proposition~\ref{sdim}. 

To prove ii), let $j = \sdim(R)$ and choose two sequences of indices
$q$ and $q'$ such that $\length(R/A_e)/q^j$ approaches $\frm$ and
$\length(R/A_{e+e'})/(qq')^j$ approaches $\frp$. Let us note that $\p +
A_{e'}^{[q]} \subset A_{e+e'}$.

Claim: $$\length (R/A_{e+e'}) \leq \length (R/\p + A_{e}^{[q']}) \leq
\length (R/A_{e}) \cdot \length (R/\p +\fm^{[q']}).$$

Indeed, the left side of the inequality is immediate since $\p
+A_{e'}^{[q]} \subset A_{e+e'}$. Set $ k = \length (R/A_{e})$ and
write a composition series $A_{e} = I_o \subset \cdots \subset I_i \subset \cdots \subset
I_k=R$. Then $(\p + A_{e}^{[q']})\subset \cdots \subset (\p+I_i)\subset \cdots   \subset
(\p+I_k)=R$.  The successive quotients of this filtration are
homomorphic images of $R/\p +\fm^{[q]}$.  Hence, $\length (R/A_{e+e'})
\le \length(R/\p + A_{e}^{[q']}) \leq k \cdot \length (R/\p
+\fm^{[q]})$.

By dividing on both sides by $(qq')^j$, with $j = \sdim(R)= \Dim(R/\p)$, and
letting $q,q'$ approach infinity we obtain the inequality stated in
ii).

For iii), note that a one-dimensional strongly $F$-regular ring $R$
is a DVR.  Part iv) follows at once from i), ii) and iii) by noting that
the Hilbert-Kunz multiplicity of a local regular ring equals one.

For v), $\Dim(R/\p) = \Dim(R)$ implies that $\p$ is a minimal prime of $R$.
Since $\p$ contains all the minimal primes of $R$, it follows that $\p$ is, in fact, the only minimal prime of $R$. But $R$ is reduced, so $\p =0$. Hence, $R$ is strongly $F$-regular.

That $\depth R \ge \sdim R$ follows from \cite{Y2}, Lemma 2.2.

\end{proof}

\begin{Q}
\label{conj}
Is is true that, for a local reduced and $F$-finite ring $\ringR$,
$\sdim(R)= \Dim(R/\p)$?
\end{Q}

Our results so far seem to indicate an affirmative answer to this
question, since we have the equality $\sdim(R) =\Dim(R/\p)$ if 
$\Dim(R/\p) \in \{ 0, 1, \Dim(R) \}$.

\begin{Prop}
Let $R=S/I$ be a Stanley-Reisner ring, then $\sum Q = \p$ where the
sum runs over all the minimal primes $Q$ of $R$.
\end{Prop}

\begin{proof}
By the definition of Stanley-Reisner rings, $I$ is a square-free
monomial ideal. In particular, $I$ can be written as an intersection
of some prime ideals in $S$, each of them generated by a string of
indeterminates.

Let us assume that the union of all the minimal primes of $R$ lifted
to $S$ involves all the indeterminates of $S$.  But then $\sum Q = \fm
\subset \p \subset \fm$ and we are done.

In the general case, let us assume that $x_1,..., x_k$ are all the
indeterminates that do not appear in any of the minimal primes of
$R$. Then one can write $R=A[x_1,...,x_k]$, with $A$ a Stanley-Reisner
ring that satisfies the condition on the minimal primes of the
previous paragraph. Hence $\fm _A = \p (A)$. It is easy to show that
$\p (R) = \p(A)R$. But then,  $\p (R) = \fm_A R$ and this equals the
sum of all the minimal primes $Q$ of $R$ as it can be easily checked.
\end{proof}

\section{The module case}

In what follows, we would like to extend some of our considerations to
the case of $R$-modules. Let $M$ be a finitely generated
$R$-module. For every $e$, we can define a new $R$-module on the
$\mathbf{Z}$-module $M$ with the $R$-multiplication given as follows:
$r \cdot m = r^{p^e}m$, for every $r \in R$ and $m \in M$. Denote this
new module by $M^{(e)}$ and, as before, let us denote $q=p^e$. The
reader should note that $M^{(e)}=M$ as abelian groups.

Denote by $a_q$ the maximal rank of a free direct summand of
$M$. Hence, $M^{(e)} = R^{a_q} \oplus M_q$ with $M_q$ an $R$-module
with no free direct summands.  Also, for every $m \in M$, denote
$\phi_{m,e}: R \to M^{(e)}$ the map that takes $r$ to $r^qm$ for every
$r \in R$. This is clearly a map of $R$-modules.

\begin{Def}
 Let $u$ be the socle generator for $\E$. For every positive integer
$e$, one can define $$A_e(M) = \{ m \in M^{(e)} : m \otimes u = 0\}$$
in $M^{(e)} \otimes \E$. Here $M^{(e)}$ is the $R$-module just
introduced above. It is easy to see that $A_e(M)$ is a submodule of
$M=M^{(e)}$ over $R$ with either of the two module structures that can be
considered.  Let us also define, $\p (M): = \cap _{e \gg 0} A_e (M)$
seen as an $R$-submodule of $M$.
\end{Def}
Note that the above Definition naturally extends  
the Definition~\ref{Aq}.

Also, $\fm^{[q]}M + \p(M) \subset A_e(M)$ for every $e$: if $r \in \fm$ and $m
 \in M$, then $r^q m  \otimes u = m \otimes ru = m \otimes 0=0$.

The following Theorem is similar to the Theorem~\ref{aq}.

\begin{Thm}

With the notations introduced above, one has that $$\length_R (M/A_e
(M))\cdot q^{\alpha(R)} = a_q.$$
\end{Thm}

\begin{proof}
Regard $A_e(M)$ as an $R$-submodule of $M$. Then $(M/A_e(M))^{(e)} =
M^{(e)}/A_e(M)$.

The proof is similar to that of Theorem~\ref{aq} and one can show that
$a_q=\length_R(M^{(e)}/A_e(M))$.  But, for every $R$-module $N$ one
has that $\length_R(N^{(e)}) = \length_R(N) \cdot q^{\alpha(R)}$.
Hence, $a_q=\length_R(M/A_e(M)) \cdot q^{\alpha(R)}.$

\end{proof}

\begin{Thm}
Let $M$ be a finitely generated $R$-module and $r\in R$ and $m \in
M$. Then $rm \in \p(M)$ if and only if $r \in \p (R)$ or $m \in \p(M)$.

\end{Thm}

\begin{proof}
The proof is similar to the that of Theorem~\ref{prime}.

Let $r \in \p(R)$ and $m \in M$. We need to show that $rm \in \p(M)$.
Consider the map $R^{(e)} \to M^{(e)}$ that sends $r'$ to $r'm$ (seen
as an element of $M$ under the usual $R$-module structure) for every
$r' \in R^{(e)}$. Tensor this map with $\E$ and let $u$ be the socle
element of this module. For every $e \gg 0$, $r \otimes u =0$.  This
maps further to zero, so $rm \otimes u=0$ for $e \gg 0$. Hence, $rm
\in \p(M)$.

For the converse implication, assume that $rm \in \p(M)$ and that $r
\notin \p(R)$ and $m \notin \p(M)$.

Since $r \notin \p(R)$, we have that for $e \gg 0$ the map
$\phi_{r,e}: R \to R^{(e)}$ splits over $R$. Similarly, for $m \notin
\p(M)$, we have that $\phi_{m,e'} : R \to M^{(e')}$ splits over $R$.

Fix $e' \gg 0$ so that $\phi_{m,e'}$ splits over $R$. Let
$q'=p^{e'}$. We will show that $r^{q'}m \notin \p(M)$ which contradicts
the fact that $rm \in \p(M)$. This will prove the converse
implication.

Take $e \gg 0$, so that $\phi_{r,e}$ splits. Define the map
$\psi_{m,e'}:R^{(e)} \to M^{(e+e')}$ by $\psi_{m,e'} (r') =
r'^{q'}m$. This is essentially the map $\phi_{m,e'}$ ``lifted'' to
$R^{(e)}$. It is an $R$-linear map and it splits over $R$. The
composition of the two maps $\phi_{r,e}$ and $\psi_{m,e'}$ defines an
$R$-linear map $R \to M^{(e+e')}$ that sends $1 \to r^{q'}m$. This map
splits over $R$ for all $e \gg 0$ (or, equivalently, $e+e' \gg 0$) and
hence $r^{q'}m \notin \p(M)$.
\end{proof}

We can define the $s$-dimension and the Frobenius splitting ratio of
an $R$-module.

\begin{Def} 
The \emph{$s$-dimension} of M, $\sdim(M)$,  is the largest integer $k$
such that $$\liminf _{q \to \infty} \frac{a_q}{q^{k+\al(R)}}$$ is not
zero. The \emph{lower Frobenius splitting ratio} of $M$, $\fr(M)$,  equals
the value of the limit introduced above, that is  $\frm(M) = \liminf
_{q \to \infty} \frac{a_q}{q^{\sdim(M)+\alpha(R)}}$. Similarly,  the
\emph{upper Frobenius splitting ratio} of $M$, $\frp(M)$,  equals $ \limsup
_{q \to \infty} \frac{a_q}{q^{\sdim(M)+\alpha(R)}}$. Whenever $\frm
(M) = \frp (M)$, we call this number the \emph{Frobenius splitting ratio} of
$M$.
\end{Def}

The following result is an extension of the case when $M=R$.

\begin{Thm}
Let $M$ be an $R$-module. Then $\sdim(M) \leq \Dim(M/\p(M))$. If $R$
is strongly $F$-regular and $M$ is an torsion free $R$-module, then
$\sdim(M) = \Dim(R) = \Dim(M)$.
\end{Thm}

\begin{proof}
From the definition of $M^{(e)}$ one can see that $(\ann_R(M^{(e)}))^q
\subset \ann_{R}(M)$, so $\Dim(M^{(e)}) = \Dim(M)$.  Write $M^{(e)} =
R^{a_q} \oplus M_q$. If $a_q \neq 0$, for some $q$, then
$\Dim(M^{(e)}) = \Dim(R)$, because $R$ injects into $M^{(e)}$.

We have seen that $a_q = \length_R (M/A_e (M))$, and $\fm^{[q]}M + \p(M) \subset A_e(M)$. So, $a_q \geq \length(M/\fm^{[q]}M + \p(M))$. Since, $e_{HK}(\fm, M/\p(M)) >0$, the first part of the Theorem follows.

If $R$ is strongly $F$-regular, then $\p(R)=0$. If $M$ is torsion-free over $R$, then $M$ injects into $R^n$ for some $n$.
Hence, if $m \otimes u = 0$ in $M^{(e)} \otimes \E$ for all $e \gg 0$, then the image of $m$ in $R^n$, say $(r_1,...,r_n)$, belongs to $\p(R)$ component-wise. However, $\p(R)=0$, and so $\p(M) =0$. So, $a_q \leq \length(M/\fm^{[q]}M)$, and since
$e_{HK} (\fm, M) >0$ we see that $\sdim(M) = \Dim (M) =\Dim(R)$.

\end{proof}

{\bf Acknowledgments}

\noindent
Much of the work for this paper has been done while the authors attended the Commutative Algebra program at the Mathematical Sciences Research Institute, Berkeley, CA, Fall of 2002. The authors would like to thank MSRI for financial support and excellent working conditions as well as the organizers of the program for making this possible.


\begin{thebibliography}{GSTT}





\bibitem{A} \textrm{I.~M.~Aberbach}, \emph{Extension of weakly and strongly $F$-regular rings by flat maps}, J. Algebra 241 (2001), 799--807.

\bibitem{AL} \textrm{I.~M.~Aberbach, G.~Leuschke}, \emph{The $F$-signature and strong $F$-regularity}, Math. Res. Lett. 10 (2003), no. 1, 51--56.

\bibitem{C} \textrm{J.~Cowden~Vassilev}, \emph{Test ideals in quotients of $F$-finite regular local rings}, Trans. Amer. Math. Soc. 350 (1998), no. 10, 4041--4051.

\bibitem{E} \textrm{F.~Enescu}, \emph{$F$-injective rings and $F$-stable primes}, Proc. Amer. Math. Soc. 131 (2003), 3379--3386.

\bibitem{F} \textrm{R.~Fedder}, \emph{F-purity and rational singularity}, Trans. Amer. Math. Soc. 278 (1983), 461--480.

\bibitem{G} \textrm{D.~Glassbrenner}, \emph{Strong $F$-regularity in images of regular rings}, Proc. Amer. Math. Soc. 124 (1996), no. 2, 345--353.

\bibitem{Ho} \textrm{M.~Hochster}, \emph{Cyclic purity versus purity in excellent Noetherian rings}, Trans. Amer. Math. Soc. 231 (1977), no. 2, 463--488.

\bibitem{HH-str} \textrm{M.~Hochster and C.~Huneke}, \emph{Tight closure and strong $F$-regularity}, Me\'moire no. 38, Soc. Math. France, 1989, 119--133.

\bibitem{HH-reg} \textrm{M.~Hochster and C.~Huneke},
  \emph{$F$-regularity, test elements, and smooth base change},
  Trans. Amer. Math. Soc. 346 (1994), no. 1, 1--62.

\bibitem{HH-adv} \textrm{M.~Hochster, C.~Huneke}, \emph{Applications of the existence of big Cohen-Macaulay algebras}, Adv. Math. 113 (1995), no. 1, 45--117. 

\bibitem{HR} \textrm{M.~Hochster, J.~L.~Roberts}, \emph{The purity of the Frobenius and local cohomology}, Adv. Math. 21 (1976), no. 2, 117--172.

\bibitem{HL} \textrm{C.~Huneke, G.~Leuschke}, \emph{Two theorems about maximal Cohen-Macaulay modules}, Math. Ann.  324  (2002),  no. 2, 391--404.

\bibitem{K} \textrm{E.~Kunz}, \emph{On Noetherian rings of characteristic $p$}, Amer. J. Math. 98 (1976), no. 4, 999--1013.

\bibitem{WY} \textrm{K.~Watanabe, K.~Yoshida}, \emph{Minimal Hilbert-Kunz multiplicity}, arXiv:math.AC/0303139, 2003.

\bibitem{Y} \textrm{Y.~Yao}, \emph{Observations on the $F$-signature of local rings of characteristic $p>0$}, preprint 2003.

\bibitem{Y2} \textrm{Y.~Yao}, \emph{Modules with finite $F$-representation
type},  preprint, 2003.

\end{thebibliography}
\end{document}